\documentclass{amsart}
\title{Rothberger bounded groups and Ramsey theory}
\author{Marion Scheepers}
\usepackage{amsfonts}
\newtheorem{theorem}{{\bf Theorem}}
\newtheorem{proposition}[theorem]{{\bf Proposition}}
\newtheorem{lemma}[theorem]{{\bf Lemma}}
\newtheorem{corollary}[theorem]{{\bf Corollary}}

\newcommand{\integers}{{\mathbb Z}}
\newcommand{\reals}{{\mathbb R}}

\newcommand{\sone}{{\sf S}_1}
\newcommand{\gone}{{\sf G}_1}

\newcommand{\open}{\mathcal{O}}
\newcommand{\onbd}{\open_{nbd}}

\address{Department of Mathematics\\ Boise State University\\ Boise, Idaho 83725}
\email{mscheepe@boisestate.edu}

\begin{document}
\maketitle
\begin{abstract} We show that
\begin{enumerate}
  \item{Rothberger bounded subgroups of $\sigma$-compact groups are characterized by Ramseyan partition relations. (Corollary 4)} 
  \item{For each uncountable cardinal $\kappa$ there is a ${\sf T}_0$ topological group of cardinality $\kappa$ such that ONE has a winning strategy in the point-open game on the group and the group is not a subspace of any $\sigma$-compact space. (Theorem \ref{rothbwonbytwo})} 
    \item{For each uncountable cardinal $\kappa$ there is a ${\sf T}_0$ topological group of cardinality $\kappa$ such that ONE has a winning strategy in the point-open game on the group and the group is $\sigma$-compact. (Corollary \ref{rothbdedtorothbomega})} 
\end{enumerate}
\end{abstract}

\section{Overview}

In topological groups boundedness properties are counterparts for covering properties in general topological spaces: Guran's notion of $\aleph_0$-boundedness is a counterpart of the Lindel\"of covering property -\cite{Guran}, while Okunev's notion of o-boundedness (named Menger-boundedness by Ko\v{c}inac, who introduced this notion independently) and Tkachenko's corresponding property of \emph{strict} o-boundedness are counterparts of $\sigma$-compactness -\cite{CH}. In this paper we consider a boundedness property which approximates Borel's metric notion of \emph{strong measure zero}. This boundedness property was introduced in unpublished work by Galvin, was later independently introduced by Ko\v{c}inac under the name of \emph{Rothberger boundedness}, and initially investigated in \cite{coc11}. 

In \cite{coc11} it was shown that a subgroup (or subset) of a metrizable topological group is Rothberger bounded if, and only if, it is strong measure zero in all left invariant metrics of the group. In \cite{coc11} we also extended some of the characterizations of strong measure zero from \cite{smzpow} to Rothberger boundedness, but we did not have techniques to also extend the Ramsey-theoretic characterization to this context. Now, in Corollary 4 of this paper, we obtain the Ramsey-theoretic characterization.

The point-open game introduced in \cite{Galvin} by Galvin is closely related to the notion of strong measure zero. Galvin proved that for Lindel\"of spaces in which each point is an intersection of countably many open sets player ONE of the point-open game has a winning strategy if, and only if, the space is countable. Apparently few, if any, uncountable examples 
of Hausdorff spaces where ONE has a winning strategy in the point-open game, have been pointed out in the literature. Another objective of this paper is to show that classical work of Comfort \cite{Comfort} and also classical work of Corson \cite{Corson} provide a wide range of topological groups that are such examples. Our analysis of these examples strengthen some results of Hernandez \cite{CH} on the theory of $\aleph_0$-bounded groups.

\section{Some terminology and notation}

For collections $\mathcal{A}$ and $\mathcal{B}$ the symbol $\sone(\mathcal{A},\mathcal{B})$ denotes the statement that
\begin{quote} For each sequence $(A_n:n<\omega)$ of elements of $\mathcal{A}$ there is a sequence $(b_n:n<\omega)$ such that for each $n$ $b_n\in A_n$, and $\{b_n:n<\omega\}\in\mathcal{B}$.
\end{quote}
Let $Y$ be a topological space. Then $\open$ denotes the collection of all open covers of $Y$. If the selection principle $\sone(\open,\open)$ holds for $Y$ we say $Y$ is a \emph{Rothberger space} (or has the \emph{Rothberger} property). This property is named after F. Rothberger who introduced it in his study of Borel's notion of strong measure zero - \cite{rothberger38}. 

A subspace $X$ of the metric space $(Y,d)$ is a \emph{strong measure zero} subspace if there is for each sequence $(\epsilon_n:n<\omega)$ of positive real numbers a partition $X=\cup_{n<\omega}X_n$ such that for each $n$ the $d$-diameter of $X_n$ is less than $\epsilon_n$. In \cite{smzpow} we characterized the strong measure zero subspaces of $\sigma$-compact metric spaces in terms of Ramseyan partition relations. 

For a topological group $(G,*)$ and a neighborhood $U$ of its identity element we define: $\open(U)=\{x*U:x\in G\}$. Then $\open(U)$ is an open cover of $G$. $\open_{nbd}$ denotes the family of all open covers of $G$ of the form $\open(U)$. A topological group is said to be \emph{Rothberger bounded} if it has the property $\sone(\open_{nbd},\open)$. For $X$ be a subset of the topological space $G$ let $\open_X$ denote the covers of $X$ by sets open in $G$. Then $X$ is said to be Rothberger bounded in $G$ if $\sone(\open_{nbd},\open_X)$ holds. In topological groups the property $\sone(\open,\open)$ is generally stronger than $\sone(\open_{nbd},\open)$. 

The symbol $\gone(\mathcal{A},\mathcal{B})$ denotes the following game of length $\omega$: Players ONE and TWO play an inning per $n<\omega$. In inning $n$ ONE first selects a member $O_n\in\mathcal{A}$, and then TWO responds by choosing a $T_n\in O_n$. A play $(O_0,\, T_0,\, \cdots,\, O_n,\, T_n,\,\cdots)$ is won by TWO if $\{T_n:n<\omega\}\in\mathcal{B}$; else, it is won by ONE.
F. Galvin \cite{Galvin} introduced the game $\gone(\open,\open)$ and proved that it is related as follows to the well-known point-open game\footnote{Since we do not need this correspondence here, we refer readers to Galvin's paper \cite{Galvin} for a definition of the point-open game.}: ONE has a winning strategy in the point-open game if, and only if, TWO has a winning strategy in $\gone(\open,\open)$. TWO has a winning strategy in the point-open game if, and only if, ONE has a winning strategy in $\gone(\open,\open)$.

\section{A Ramseyan characterization of some Rothberger bounded groups.}

Let $\mathcal{A}$ and $\mathcal{B}$ be families of sets and let $n$ and $k$ be positive integers. The symbol 
\[
  \mathcal{A}\rightarrow(\mathcal{B})^n_k
\]
denotes the statement that for each $A\in\mathcal{A}$ and for each function $f:\lbrack A\rbrack^n\rightarrow\{1,\,\cdots,\, k\}$ there is a $B\subseteq A$ and an $i\in\{1,\,\cdots,\,k\}$ such that $B\in\mathcal{B}$ and $f$ is constant of value $i$ on $\lbrack B\rbrack^n$.

An open cover $\mathcal{U}$ of $Y$ is said to be an $\omega$-cover if $Y\not\in\mathcal{U}$, but for each finite set $F\subseteq Y$ there is a $U\in\mathcal{U}$ with $F\subseteq U$. $\Omega$ denotes the family of all open $\omega$-covers of $Y$.

Theorem 9 of \cite{smzpow} gives the following characterization of strong measure zero subsets of $\sigma$-compact metric spaces:
\begin{theorem}\label{smzpowTh9} For $X$ a subspace of a $\sigma$-compact metric space $Y$ the following are equivalent:
\begin{enumerate}
  \item{$Y$ has the property $\sone(\open,\mathcal{O}_X)$.} 
  \item{$X$ has strong measure zero (in all equivalent metrics on $Y$).}
  \item{ONE has no winning strategy in the game $\gone(\open,\mathcal{O}_X)$.}  
  \item{For each positive integer $k$, $\Omega\rightarrow(\open_X)^2_k$.} 
\end{enumerate}
\end{theorem}

Failure of Lebesgue's Covering Lemma (which holds for compact metrizable spaces) was the main obstacle towards extending Theorem \ref{smzpowTh9} beyond metrizable $\sigma$-compact spaces. We found two non-metric situations in which an appropriate generalization of Lebesgue's Covering Lemma holds. Here is the first of the two\footnote{Lemma \ref{lebesguecov} must be well-known but I was not able to track down a reference for it.}: 
\begin{lemma}\label{lebesguecov} Let $(G,*)$ be a ${\sf T}_0$ topological group, let $\mathcal{U}$ be an open cover of $G$ and let $K$ be a nonempty compact subset of $G$. Then there is a neighborhood $N$ of the identity of $G$ such that for each $x\in K$ there is a $U\in\mathcal{U}$ with $x*N\subseteq U$.
\end{lemma}
{\flushleft{\bf Proof:}} Let $\mathcal{U}$ be an open cover of $G$. For each $x\in K$ choose a $U(x)\in\mathcal{U}$ with $x\in U(x)$. Then choose for each $x$ a neighborhood $N_x$ of $G$'s identity such that $x*{N_x}^2\subseteq U(x)$. Now $\{x*N_x:x\in K\}$ is an open cover of $K$. Since $K$ is compact this cover has a finite subset that covers $K$, say $\{x_1*N_{x_1},\cdots,x_k*N_{x_k}\}$. Define $N$ by 
\[
  N=N_{x_1}\cap\cdots\cap N_{x_k}.
\]
Then $N$ is as required. For consider any $k\in K$. Choose $i$ so that $k\in x_i*N_{x_i}$. Now consider any $q\in k*N$. Fix $a\in N$ with $q = k*a$. Then we have $q\in x_i*N_{x_i}*a \subseteq x_i*N^2_{x_i} \subseteq U(x_i)$. Since $q$ was an arbitrary element of $k*N$ we find that $k*N\subseteq U(x_i)$. $\Box$

Here is the version of Theorem \ref{smzpowTh9} for $\sigma$-compact topological groups:
\begin{theorem}\label{rboundedinsigmacompact}
Let $(G,*)$ be a $\sigma$-compact topological group and let $X$ be a subset of $G$. The following are equivalent:
\begin{enumerate}
  \item{$(G,*)$ satisfies $\sone(\open,\open_X)$.} 
  \item{$(G,*)$ satisfies $\sone(\open_{nbd},\open_X)$.} 
  \item{$X$ is a Rothberger bounded subset of $(G,*)$.} 
  \item{ONE has no winning strategy in the game $\gone(\open,\open_X)$.} 
  \item{For each positive integer $k$, $\Omega\rightarrow (\open_X)^2_k$} 
\end{enumerate}
\end{theorem}
{\flushleft{\bf Proof:}} The implications $(2)\Rightarrow (3)$, $(4)\Rightarrow (5)$ and $(5)\Rightarrow(1)$ follow the ideas of the proof of Theorem 9 of \cite{smzpow}, with some standard modifications. We prove $(3)\Rightarrow(4)$ here since its proof uses the perhaps new Lemma \ref{lebesguecov}.

{\flushleft{\underline{\bf $(3)\Rightarrow(4):$}} } Let $F$ be a strategy for player ONE in the game $\gone(\open,\open_X)$ on $G$. Since $G$ is $\sigma$-compact, write $G=\bigcup_{n<\omega}G_n$ where for each $n$ we have $id_G\in G_n\subseteq G_{n+1}$ and $G_n$ is compact. For each $n$, $X_n = X\cap G_n$ is a Rothberger bounded subset of $G$. To defeat ONE's strategy TWO will concentrate attention on specific $X_n$'s in specific innings. To this end, partition $\omega$ into infinitely many infinite subsets $S_n$. For innings numbered by members of $S_n$ TWO will focus on $X_n$.

We now use ONE's strategy $F$ to recursively define a sequence $(N_k:k<\omega)$ and an array $(\mathcal{U}(T_0,\cdots,T_k):k<\omega)$ where
\begin{enumerate}
  \item{For each $k$, $N_k$ is a \emph{symmetric}\footnote{\emph{i.e.}, $N_k^{-1}=N_k$.} neighborhood of the identity of $G$;}
  \item{With $n_0$ such that $0\in S_{n_0}$, $\mathcal{U}(\emptyset)$ is a finite subset of $F(\emptyset)$ (ONE's first move) which covers $G_0$, and $N_0$ is a neighborhood of the identity of $G$ such that for each $x\in G_0$ there is a $V\in\mathcal{U}(\emptyset)$ with $x*(N_0*N_0)\subseteq V$.}
  \item{For each $(T_0,\cdots,T_k)$ such that $T_0\in F(\emptyset)$, $T_1\in F(T_0)$, $\cdots$ and $T_k\in F(T_0,\cdots,T_{k-1})$ and for $n_{k+1}$ such that $k+1\in S_{n_{k+1}}$ we have $\mathcal{U}(T_0,\cdots,T_k)$ a finite subset of $F(T_0,\cdots,T_k)$ that covers $G_{n_{k+1}}$. Note that there are only finitely many such $(T_0,\cdots,T_k)$. $N_{k+1}$ is a neighborhood of the identity of $G$ such that for each such sequence $(T_0,\cdots,T_k)$ and for each $x\in G_{n_{k+1}}$ there is a $U\in\mathcal{U}(T_0,\cdots,T_{k})$ with $x*N_{k+1}*N_{k+1}\subseteq U$.}
\end{enumerate}

With this data available, construct a play against $F$ won by TWO as follows: Fix an $m<\omega$. Since $X_m$ is Rothberger bounded select for each $k\in S_m$ an $x_k\in G$ such that $(x_k*N_k:k\in S_m)$ covers $X_m$. 

We may assume each $x_k$ is in $X_m$ - for suppose an $x_k$ is not in $X_m$. If $x_k*N_k \cap X_k = \emptyset$, we may with impunity replace this $x_k$ by one from $X_k$. However, if $x_k*N_k\cap X_k\neq\emptyset$, then let $y$ be an element of this intersection. We claim that $x_k*N_k\cap X_k\subseteq y*N_k^2$. For let $z\in x_k*N_k\cap X_k$ be given. Choose $v$ in $N_k$ with $z = x_k*v$. Since we also have $y\in x_k*N_k$, choose $w\in N_k$ with $y=x_k*w$. Now $x_k = y*w^{-1} \in N_k*N_k$ and so $z = y*(w^{-1}*v)\in y*N_k^2$.

Thus, we may choose for each $k\in S_m$ an $x_k\in X_k$ such that $(x_k*U_k^2:k\in S_m)$ covers $X_k$.
Finally, recursively choose a sequence $(T_k:k<\omega)$ as follows: Choose $T_0\in \mathcal{U}(\emptyset)$ with $x_0*N_0^2 \subseteq T_0$. With $T_0,\, \cdots,\, T_m$ chosen, choose $T_{m+1}\in\mathcal{U}(T_0,\cdots,T_m)$ with $x_{m+1}*N_{m+1}^2\subseteq T_{m+1}$. Then the sequence
\[
  F(\emptyset),\, T_0,\, F(T_0),\, \cdots,\, T_k,\, F(T_0,\cdots,T_k),\, T_{k+1},\, \cdots
\]
is an $F$-play lost by ONE.

In the proof of (5) implies (1) we use that $\sone(\open,\open_X)$ is equivalent to $\sone(\Omega,\open_X)$.
$\Box$

Note that $(4)$ also implies that ONE has no winning strategy in the game $\gone(\open_{nbd},\open_X)$ on $(G,*)$, which in turn implies $(2)$. In particular we have the following Corollary:
\begin{corollary}\label{rothbbdedgp}
Let $(H,*)$ be a subgroup of a $\sigma$-compact group $(G,*)$. The following are equivalent:
\begin{enumerate}
  \item{$(G,*)$ satisfies $\sone(\open,\open_H)$.} 
  \item{$(H,*)$ is Rothberger bounded.} 
  \item{ONE has no winning strategy in the game $\gone(\open,\open_H)$.}
  \item{On $H$ ONE has no winning strategy in the game $\gone(\open_{nbd},\open)$.}  
  \item{For each positive integer $k$, $(G,*)$ satisfies $\Omega\rightarrow (\open_H)^2_k$} 
\end{enumerate}
\end{corollary}
{\flushleft{\bf Proof:}} (1) implies that $(G,*)$ satisfies $\sone(\open_{nbd},\open_H)$. But then by Theorem 13 of \cite{coc11}, $(H,*)$ satisfies $\sone(\open_{nbd},\open)$ in the relative topology, proving (2). Since (2) states that $H$ is a Rothberger bounded subset of $G$, Theorem \ref{rboundedinsigmacompact} gives the implication from (2) to (3) and from (3) to (5). Also, (3) implies (4) which implies (2). The proof that (5) implies (1) uses the fact that $\sone(\open,\open_H)$ is equivalent to $\sone(\Omega,\open_H)$.
$\Box$\\

Corollary \ref{rothbbdedgp} improves Theorem 22 of \cite{coc11} in that it does not require the group $(G,*)$ to be metrizable. One might wonder how widely applicable Corollary \ref{rothbbdedgp} really is. We shall address this in the next two sections by showing that:
\begin{itemize}
  \item{For each infinite cardinal number $\kappa$ there is a ${\sf T}_0$ Rothberger bounded group $(G,*)$ of cardinality $\kappa$ which is not a subgroup of any $\sigma$-compact group, and TWO has a winning strategy in the game $\gone(\open_{nbd},\open)$ on $(G,*)$.} 
  \item{For each infinite cardinal number $\kappa$ there is a ${\sf T}_0$ $\sigma$-compact Rothberger bounded group of cardinality $\kappa$ for which TWO has a winning strategy in the game $\gone(\open,\open)$.}
\end{itemize}

\section{Rothberger groups not embedding into any $\sigma$-compact space.}

Not every ${\sf T}_0$ Rothberger bounded group is a subgroup of a $\sigma$-compact group, as illustrated by the following example of Comfort and Ross (\cite{CR}, Example 3.2). We precede the example with a few general remarks about $P$-spaces. A topological space is said to be a $P$-space if each ${\sf G}_{\delta}$ set is open. Evidently, every subspace of a $P$-space is a $P$-space. Every countably infinite subspace of a ${\sf T}_2$ $P$-space is closed and discrete. It follows that a compact $P$-space is finite, and thus a $\sigma$-compact $P$-space is countable. Thus, no uncountable Lindel\"of $P$-space is a subspace of a $\sigma$-compact ${\sf T}_2$-space. If a topological group $(G,*)$ is a Lindel\"of $P$-space then it is Rothberger bounded in a strong sense: Let $(U_n:n<\omega)$ be a sequence of neighborhoods for the identity. Then $U = \cap_{n<\omega}U_n$ is a neighborhood for the identity. Since the group is $\aleph_0$ bounded fix a sequence $(x_n:n<\omega)$ of elements of the group such that $x_n*U,\, n<\omega$ covers the group. Then the sequence $(x_n*U_n:n<\omega)$ witnesses that the group is Rothberger bounded. 

We now define the example: The underlying set of the group $G$ is
\[
  G:=\{f\in\,^{\omega_1}2:\vert\{\alpha: f(\alpha)\neq 0\}\vert<\aleph_0\};
\]
Endow $G$ with the {\sf G}$_{\delta}$ topology. \cite{CR} shows that $(G,\oplus)$ is a Lindel\"of $P$-group (and thus ${\sf T}_4$). Thus $(G,\oplus)$ is an uncountable Rothberger bounded group that is not contained in a $\sigma$-compact group.
Theorem 2.3 of Comfort in \cite{Comfort} implies the following generalization of the above example:
\begin{theorem}[Comfort]\label{Comfortcountable}
Let $(G_i,*_i)$, $i\in I$, be a family  of countable topological groups. Endow the product $\prod_{i\in I}G_i$ with the countable box topology. Then the subgroup
\[
  G:=\{f\in\prod_{i\in I}G_i: \vert\{j\in I:f(j)\neq {\sf id}_j\}\vert<\aleph_0\}
\]
is a Lindel\"of $P$-group.
\end{theorem}
In particular, for each uncountable cardinal number $\kappa$ there is a ${\sf T}_0$ Lindel\"of $P$ group of cardinality $\kappa$.
One can prove an analogue of Theorem \ref{rboundedinsigmacompact} also for Lindel\"of  
$P$-groups, because: Galvin proved a result that implies that if a space is a Lindel\"of P-space then it is a Rothberger space - see the Lemma in Section 2 of \cite{GN}. Thus, Lindel\"of $P$ groups satisfy the stronger selection principle $\sone(\open,\open)$. But the following are equivalent:
\begin{enumerate}
  \item{Topological space $X$ satisfies $\sone(\open,\open)$;}
  \item{ONE has no winning strategy in the game $\gone(\open,\open)$ played on $X$;}
  \item{For each positive integer $k$, $\Omega\rightarrow(\open)^2_k$ holds for $X$.}   
\end{enumerate}
The equivalence of (1) and (2) is due to Pawlikowski \cite{JP}, and the equivalence with (3) was proved in \cite{msrothbramsey}.
Using ideas of \cite{CH} we now show:
\begin{proposition}\label{twowinsrothberger} For each uncountable cardinal $\kappa$ there is a ${\sf T}_0$ Rothberger bounded group of cardinality $\kappa$ such that TWO has a winning strategy in $\gone(\open_{nbd},\open)$.
\end{proposition}

In the proof we will make use of the following elementary game of $\omega$ innings played on an infinite set $S$: In inning $n$ player ONE chooses a countable subset $W_n$ of $S$ and TWO responds by choosing a point $b_n\in  W_n$. ONE must further obey the rule that for each $n$, $W_n\subseteq W_{n+1}$. A play $(W_0,\, b_0,\, \cdots,\, W_n,\, b_n,\,\cdots)$ is won by TWO if for each $x\in\bigcup_{n<\omega}W_n$ there are infinitely may $n$ with $b_n = x$. A standard argument shows that TWO has a winning perfect information strategy in this game. Call this game the ``countable - one game".

{\flushleft{\bf Proof of Proposition \ref{twowinsrothberger}:}} Let $\kappa$ be an uncountable cardinal number. Let $(G_{\alpha}:\alpha<\kappa)$ be a sequence of discrete  countable groups and define $G$ to be the direct product 
\[
  \prod_{\alpha<\kappa}G_{\alpha}
\]
endowed with the countable box topology. 

Then the subset
  $G^* = \{f\in G:\vert \{\alpha: f(\alpha)\neq {\sf id}_{\alpha}\}\vert<\aleph_0\}$
endowed with the relative topology is by Theorem \ref{Comfortcountable} a Lindel\"of ${\sf P}$-group.
For a countable set $B\subset\kappa$, let $\Pi_B$ denote the projection of $G$ onto $\prod_{\alpha\in B}G_{\alpha}$. Then the  set
\[
  U_B = G^* \bigcap \Pi^{\leftarrow}_B\lbrack\{{\sf id}_B\}\rbrack 
\]
is a basic neighborhood of the identity element of $G^*$. Also
  $\mathcal{D} = \{U_B:B\in\lbrack\kappa\rbrack^{\leq\aleph_0}\}$
is a neighborhood basis for the identity element of $G^*$, and each $U_B$ is a subgroup of the group $G^*$. Since $G^*$ is a Lindel\"of ${\sf P}$-group, each open cover of the form $\mathcal{O}(U_B)$ has a countable subcover, and this means that the subgroup $U_B$ has countably many distinct left cosets in $G^*$. 

Now we show that TWO has a winning strategy in the game $\gone(\onbd,\open)$ played on $G^*$. Since $\mathcal{D}$ is a neighborhood basis of the identity element of $G^*$, we may assume that for each $n<\omega$ ONE's $n$-th move is of the form $\open(U_{B_n})$, $B_n$ a countable subset of $\kappa$. And since TWO may replace ONE's move $\open(U_B)$ with a move $\open(U_C)\subset \open(U_B)$ and respond to the replacement move instead, we may further assume that ONE's moves are such that for each $n$, $B_n\subseteq B_{n+1}$, that is, $U_{B_{n+1}}\subseteq U_{B_n}$.

Also, for each move $\open(U_{B_n})$ by ONE, TWO chooses a countable set $A_n\subseteq G^*$ such that $\{x* U_{B_n}: x\in A_n\}$ is the set of distinct left cosets of $U_{B_n}$ in $G^*$. Since TWO has perfect information and for each $n$ $U_{B_{n+1}}\subseteq U_{B_n}$, TWO may select the sets $A_n$ such that for each $n$ we have $A_n\subseteq A_{n+1}$. We may assume for each $n$ that for each $x\in A_n$, if $x(\alpha)\neq {\sf id}_{G_{\alpha}}$ then $\alpha\in B_n$.

Now let $F$ be a winning perfect information strategy for TWO in the countable-1 game on $G^*$. We define a strategy $\sigma$ for TWO in $\gone(\onbd,\open)$ as follows:

Given ONE's move $\open(U_{B_0})$ in $\gone(\onbd,\open)$ TWO first fixes $A_0$ as above, considered as a move of ONE of the countable-1 game. Then in that game TWO moves $x_0 = F(A_0)$. Then TWO responds in $\gone(\onbd,\open)$ with $\sigma(\open(U_{B_0}) = x_0*U_{B_0}$. 

In the next inning of $\gone(\onbd,\open)$ ONE moves $\open(U_{B_1})$. TWO first fixes the countble set $A_1$ as above and consider it as a move of ONE in the countable - 1 game on $G^*$. In that game TWO moves $x_1 = F(A_0,A_1)$. Then TWO responds in $\gone(\onbd,\open)$ with $\sigma(\open(U_{B_0}),\open(U_{B_1})) = x_1*U_{B_1}$, and so on, as depicted in the following diagram.

\begin{center}
\begin{tabular}{c|l}
\multicolumn{2}{c}{$\gone(\onbd,\open)$}\\ \hline\hline
\multicolumn{2}{c}{}\\
ONE              & TWO                                     \\ \hline
$\open(U_{B_0})$   &                                       \\
                 &                                         \\
                 & $x_0*U_{B_0}$ \\
$\open(U_{B_1})$   &                                       \\
                 &                                         \\
                 & $x_1*U_{B_1} $ \\
$\vdots$         & $\vdots$ \\
\end{tabular}
\hspace{1in}
\begin{tabular}{c|l}
\multicolumn{2}{c}{Countable-1 game}                        \\ \hline\hline
\multicolumn{2}{c}{}\\ 
ONE   & TWO                                                 \\ \hline
      &                                                     \\
$A_0$ & $x_0 = F(A_0)$                                      \\
      &                                                     \\
      &                                                     \\
$A_1$ & $x_1=F(A_0,A_1)$                                    \\
      &                                                     \\
$\vdots$         & $\vdots$ \\
\end{tabular}
\end{center}

We shall now see that $\sigma$ is a winning strategy for TWO in the game $\gone(\onbd,\open)$. Thus, consider a $\sigma$-play
\[
  \open(U_{B_0}),\, x_0*U_{B_0},\,\open(U_{B_1}),\, x_1*U_{B_1},\,\cdots,\, \open(U_{B_n}),\, x_n*U_{B_n},\,\cdots
\]

Let $x\in G^*$ be given. We must show that 
  $x\in\bigcup_{n<\omega}x_n*U_{B_n}$.
Put $B=\cup_{n<\omega}B_n$. Then evidently $U_B\subseteq U_{B_n}$ holds for each $n$. Also
  $L = x*U_B$
is  left coset of $U_B$ in $G^*$.

{\flushleft{\bf Claim 1:}} For each $n<\omega$ there is a $y\in A_n$ with $L\subseteq y*U_{B_n}$.

To see this, fix $n<\omega$. Since $U_B\subseteq U_{B_n}$ we have $x*U_B\subseteq x*U_{B_n}$. But $x*U_{B_n}$ is a left coset of $U_{B_n}$ in $G^*$, and so by the choice of $A_n$ there is a $y\in A_n$ with $x*U_{B_n} = y*U_{B_n}$.

Thus, choose for each $n$ a $y_n\in A_n$ such that $L\subseteq y_n*U_{B_n}$.

{\flushleft{\bf Claim 2:}} For each $n<\omega$ we have $y_{n+1}*U_{B_{n+1}}\subseteq y_n*U_{B_n}$.

For suppose on the contrary that $y_{n+1}*U_{B_{n+1}}\not\subseteq y_n*U_{B_n}$. Then we also have $y_{n+1}*U_{B_n}\not\subseteq y_n*U_{B_n}$, so that these are distinct left cosets of $U_{B_n}$ in $G^*$, and so are disjoint. But this contradicts the fact that $\emptyset\neq L\subseteq y_{n+1}*U_{B_{n+1}}\cap y_n*U_{B_{n}}$.

Towards proving the next claim first note that if $x*U_{B_n}\subseteq y*U_{B_n}$ then $(\forall \alpha\in B_n)(x(\alpha)=y(\alpha))$.

{\flushleft{\bf Claim 3:}} For each $n<\omega$ we have $support(y_n)\subseteq support(y_{n+1})$.

For by Claim 2 we find a $u\in U_{B_n}$ such that $y_{n+1} = y_n*u$. Now $support(y_n)\cap support(u) = \emptyset$, and so $support(y_{n+1}) = support(y_n)\cup support(u)\supseteq support(y_n)$.

{\flushleft{\bf Claim 4:}} For each $n<\omega$ we have $support(y_n)\subseteq support(x)$.

Since $x\in L\subseteq y_n*U_{B_n}$ it follows that $y_n^{-1}*x\in U_{B_n}$, and so for each $\alpha\in B_n$ we have $x(\alpha) = y_n(\alpha)$. Since $support(y_n)\subseteq B_n$, Claim 4 follows.

Since $x$ is in $G^*$ it has finite support. Claims 3 and 4 imply that there is a $k$, from now on fixed, such that for all $n\geq k$ we have $support(y_n) = support(y_k)$. It follows that for all $n\ge k$, $y_n = y_k$. But then, for all $n\ge k$, $x\in y_k*U_{B_n}$, which implies that $x\in y_k*U_B$. But for infinitely many $n$ we have $x_n = y_k$,a nd so for such an $n$ larger than $k$, $x\in x_n*U_{B_n}$.

This completes the proof that $\sigma$ is a winning strategy for TWO.
$\Box$

From this we now derive that TWO in fact has a winning strategy in the game $\gone(\open,\open)$, typically a harder game for player TWO. Towards this we need another generalization of the Lebesgue Covering Lemma, this time for Lindel\"of P-groups:
\begin{proposition}\label{pgpLebesgue} Let $(G,*)$ be a Lindel\"of P-group. Then there is for each open cover $\mathcal{U}$ of $G$ a neighborhood $N$ of the identity of $G$ such that for each $x\in G$ there is a $U\in\mathcal{U}$ such that $x*N\subseteq U$. 
\end{proposition}
{\flushleft{\bf Proof:}} Let $\mathcal{U}$ be an open cover of $G$. For each $x\in G$ choose a neighborhood $U_x$ of $G$'s identity such that $U_x$ is an open subgroup of $G$, and there is a $U\in\mathcal{U}$ with $x*U_x\subseteq U$. Since $G$ is Lindel\"of we find $x_n,\, n<\omega$ such that $\mathcal{F} = \{x_n*U_{x_n}:n<\omega\}$ is an open cover of $G$ and refines $\mathcal{U}$. Now since $G$ is a $P$-space, choose an open neighborhood $N$ of the identity such that $N$ is a subgroup of $G$, and $N\subseteq \bigcap_{n<\omega}U_{x_n}$. 

We claim that $N$ is as required. For consider any $x\in G$. Then $x*N$ is a left coset of $N$ in $G$. We claim there is an $n$ with $x*N\subseteq x_n*U_{x_n}$. For if not, then for each $n$ we have $x*N\not\subseteq x_n*U_{x_n}$. But we have $x*N\subseteq x*U_{x_n}$, a left coset of $U_{x_n}$ in the group $G$. Thus $x*U_{x_n}\neq x_n*U_{x_n}$, and as $x_n*U_{x_n}$ is also a left coset of $U_{x_n}$, we have $x*N\cap x_n*U_{x_n} = \emptyset$. But then the family $\{x_n*U_{x_n}:n<\omega\}$ does not cover the subset $x*N$ of $G$, contradicting the fact that $\mathcal{F}$ is a cover of $G$.
$\Box$

\begin{theorem}\label{rothbwonbytwo}
For each infinite cardinal $\kappa$ there is a ${\sf T}_0$ Lindel\"of $P$-group of cardinality $\kappa$ such that TWO has a winning strategy in the game $\gone(\open,\open)$.
\end{theorem}
{\flushleft{\bf Proof:}} Let $(G,*)$ be a Lindel\"of $P$-group for which TWO has a winning strategy in $\gone(\open_{nbd},\open)$ (as for example in Proposition \ref{twowinsrothberger}). Let $F$ be TWO's winning strategy in that game.

Define a strategy $\sigma$ for TWO in the game $\gone(\open,\open)$ as follows: When ONE plays the open cover $\mathcal{U}_1$, choose a neighborhood $N_1$ of the identity of $G$ as in Proposition \ref{pgpLebesgue}, and then let $\sigma(\mathcal{U}_1)$ be an element $U$ of $\mathcal{U}_1$ such that $F(N_1)*N_1\subseteq U$. When ONE plays the next open cover $\mathcal{U}_2$ choose a neighborhoood $N_2$ of the identity  of $G$ as in Proposition \ref{pgpLebesgue}, and then let $\sigma(\mathcal{U}_1,\mathcal{U}_2)$ be an element $U$ of $\mathcal{U}_2$ such that $F(N_1,N_2)*N_2\subseteq U$, and so forth. 

Then $\sigma$ is a winning strategy for TWO.
  $\Box$

Call an open cover $\mathcal{U}$ of a topological space a $\gamma$-cover if $\mathcal{U}$ is infinite, and each infinite subset of $\mathcal{U}$ still covers the space. The symbol $\Gamma$ denotes the collection of open $\gamma$ covers of a space. In \cite{GN} Gerlits and Nagy introduced the notion of a $\gamma$-space: A topological space which satisfies the selection principle $\sone(\Omega,\Gamma)$ is said to be a $\gamma$-\emph{space}. It is evident that each $\gamma$-space is a Rothberger space. In Theorem 1 of \cite{GN} the authors prove
\begin{theorem}[Gerlits-Nagy] For a ${\sf T}_{3\frac{1}{2}}$-space TWO has a winning strategy in $\gone(\open,\open)$ if, and only if, TWO has a winning strategy in $\gone(\Omega,\Gamma)$.
\end{theorem}
It is also evident that if TWO has a winning strategy in the game $\gone(\Omega,\Gamma)$, then the underlying space is a $\gamma$-space.
\begin{corollary}\label{GNcorollary} For each uncountable cardinal number $\kappa$ there is a Lindel\"of $P$-group of cardinality $\kappa$ on which TWO has a winning strategy in the game $\gone(\Omega,\Gamma)$.
\end{corollary}
{\flushleft{\bf Proof:}} By a result of Gerlits and Nagy if TWO has a winning strategy in $\gone(\open,\open)$, then TWO has a winning strategy in $\gone(\Omega,\Gamma)$.
$\Box$

\section{Large $\sigma$-compact Rothberger bounded ${\sf T}_0$ groups.}

In Proposition 4 of \cite{Corson}, Corson proves essentially the following theorem\footnote{Corson formulates the proposition for the case when the factor spaces $X_i$ are all the real line. But the argument gives the more general result of Theorem \ref{Corson}.}:
\begin{theorem}[Corson]\label{Corson} Let $\{X_i:i\in I\}$ be a family of $\sigma$-compact topological groups and for each $i$ let $e_i$ be the identity element of $X_i$. Then the subgroup 
\[
  G:=\{f\in\prod_{i\in I}X_i:\vert\{j\in I:f(j)\neq e_j\}\vert< \aleph_0\}
\]
is $\sigma$-compact.
\end{theorem}

\begin{corollary}\label{corsoncor} For each infinite cardinal number $\kappa$ there is a ${\sf T}_0$ $\sigma$-compact Rothberger bounded group of cardinality $\kappa$. 
\end{corollary}
{\flushleft{\bf Proof:}} Let cardinal number $\kappa$ be given, and take $I$ to be $\kappa$. For each $i\in I$ take $X_i$ to be $\integers$, the additive group of integers. Now consider the group $G$ as in Corson's Theorem. $G$ is in fact Rothberger bounded. To see this let for each $n$ a neighborhood $U_n$ of the identity element of $G$ be given. We may assume each $U_n$ is a basic open set, and thus that there is a finite set $F_n\subseteq I$ and for each $i\in F_n$ a neighborhood $N_{n,i}$ of $e_i$ such that $U_n$ is of the form 
  $\{f\in G:(\forall i\in F_n)(f(i)\in N_{n,i})\}$.
Now $C = \cup_{n<\omega}F_n$ is a countable subset of $I$ and $G_C = \{f\lceil_C:f \in G\}$ is evidently Rothberger bounded in $\prod_{i\in C}X_i$. For each $n$ choose a $g_n\in G_C$ such that $G_C\subseteq \cup_{n<\omega}g_n*U_n\lceil_C$. For each $n$ choose $f_n\in G$ with $f_n\lceil_C = g_n$. Then it follows that $G\subseteq \cup_{n<\omega}f_n*U_n$. Thus $G$ is a $\sigma$-compact Rothberger bounded group of cardinality $\kappa$. $\Box$

Using the method of proof of Proposition \ref{twowinsrothberger}, one proves
\begin{proposition}\label{twowinssigmacompactrothberger}
In the groups of Corollary \ref{corsoncor} TWO has a winning strategy in the game $\gone(\onbd,\open)$.
\end{proposition}
{\flushleft{\bf Proof}} Note that for each finite subset $B$ of $\kappa$ the set $U_B =\{f\in G:(\forall \alpha\in B)(f(\alpha) = id_{G_{\alpha}}\}$ is in fact a subgroup of $G$, is a neighborhood of the identity element of $G$, and the set of such $U_B$ form a neighborhood basis of the identity element of $G$. Now apply the argument of Proposition \ref{twowinsrothberger}. $\Box$

In fact the groups of Corollary \ref{corsoncor} satisfy the stronger selection principle $\sone(\open,\open)$:
One can prove more, namely
\begin{theorem}\label{twowinssigmacompact} In the groups of Corollary \ref{corsoncor} TWO has a winning strategy in the game $\gone(\open,\open)$ on $G$.
\end{theorem}
{\flushleft{\bf Proof: }} Let $(G,*)$ be such a group and $F$ be a winning strategy for TWO in the game $\gone(\onbd,\open)$ on this group. Since the group is $\sigma$-compact we  write $G = \bigcup_{n=0}^{\infty}G_n$ where for each $n$ $G_n\subseteq G_{n+1}$ and $G_n$ is compact.

For each open cover $\mathcal{U}$ of $G$ and for each $n$ choose a neighborhood $U(\mathcal{U},n)$ of the identity element of $G$ such that $U(\mathcal{U},n)$ is a subgroup of $G$ and for each $x\in G_n$ there is a $V\in\mathcal{U}$ such that $x* U(\mathcal{U},n) \subseteq V$. Let $V(x,\mathcal{U},n)$ be such a $V$.

Define a strategy $\sigma$ for TWO of $\gone(\open,\open)$ as follows: When ONE plays the open cover $O_1$ in the first inning, TWO simulates a move for ONE in $\gone(\onbd,\open)$ as $\open(U(O_1,1))$, applies the winning strategy $F$ to this move to obtain $x_1*U(O_1,1) = F(\open(U(O_1,1)))$. Then if $G_1\cap x_1*U(O_1,1)\neq \emptyset$ we find for an $x\in G_1$ that $x*U(O_1,1) = x_1*U(O_1,1)$, a left coset of $U(O_1,1)$ in $G$. Then TWO fixes such an $x$ and plays
\[
  \sigma(O_1) = V(x,O_1,1) \in O_1.
\]
If $G_1\cap x_1*U(O_1,1) = \emptyset$ then TWO chooses an arbitrary element $x\in G_1$ and plays 
\[
  \sigma(O_1) = V(x,O_1,1) \in O_1.
\]
When ONE next moves $O_2$, TWO simulates a move for ONE in $\gone(\onbd,\open)$ as follows: Define $\open = \{U\cap V:U\in O_1 \mbox{ and }V\in O_1\}\setminus\{\emptyset\}$ and let ONE's move be $\open(U(\open,2))$. TWO's response using $F$ is $x_2*U(\open,2)=F(\open(U(O_1,1)),\open(U(\open,2))$. Consider $x_2*U(\open,2)\cap G_2$. If this is nonempty select any $x$ in this intersection. as before we have $x*U(\open,2)=x_2*U(\open,2)$, and now TWO responds with 
\[
  \sigma(O_1,O_2) = V(x,\open,2)\in O_2.
\]
If on the other hand the intersection is empty then TWO chooses any $x\in G_2$ and responds with
\[
  \sigma(O_1,O_2) = V(x,\open,2)\in O_2.
\]
This procedure describes a strategy for TWO in the game $\gone(\open,\open)$ on $G$.

To see that $\sigma$ is a winning strategy for TWO, consider any $\sigma$-play
\[
  O_1,\, \sigma(O_1),\, O_2,\,\sigma(O_1,O_2),\,\cdots
\]
Let an $x\in G$ be given. Fix the least $m$ with $x\in G_m$. By the definition of $\sigma$ we have an associated sequence 
\[
  U_n = U(\open_n,n)
\]
of subgroups of $G$ that are neighborhoods for the identity element where for each $n$ we have $U_{n+1}\subset U_n$, and $\open_1 = O_1$ while $\open_{n+1} = \{U\cap V:\, U\in \open_n \mbox{ and }V\in O_{n+1}\}\setminus\{\emptyset\}$, and elements $x_n$ of $G$ such that 
\[
  x_1*U_1 = F(\open(U_1)) \mbox{ and } x_{n+1}*U_{n+1} = F(\open(U_1),\cdots,\open(U_n)).
\]

But then this is an $F$-play of $\gone(\onbd,\open)$ and thus won by TWO, meaning there are infinitely many $n$ with $x\in x_n*U_n$. Thus, fix an $n>m$ with $x\in x_n*U_n$. But then $G_n\cap x_n*G_n\neq \emptyset$ and thus $\sigma(O_1,\cdots,O_n)\supseteq x*U_n$, meaning $x\in \sigma(O_1,\cdots,O_n)$ $\Box$

\begin{corollary}\label{sigmacompactGN} For each infinite cardinal number $\kappa$ there is a $\sigma$-compact ${\sf T}_0$ topological group of cardinality $\kappa$ such that TWO has a winning strategy in the game $\gone(\Omega,\Gamma)$.
\end{corollary}

\begin{proposition}\label{rothbdedtorothb} Let $(G,*)$ be a $\sigma$-compact ${\sf T}_0$ topological group with property $\sone(\onbd,\open)$. Then $G$ has the property $\sone(\Omega,\Gamma)$.
\end{proposition} 
{\flushleft{\bf Proof:}} 
This follows directly from Corollary \ref{rothbbdedgp} part (1) that $(G,*)$ has property $\sone(\open,\open)$ since in the notation of that corollary $G=H$ and $\open_H=\open$. 

Recall that a ${\sf T}_0$ topological group is ${\sf T}_{3\frac{1}{2}}$. Since a compact ${\sf T}_3$-space is Rothberger if, and only if, it is scattered (\cite{BCM} Proposition 34), if, and only if, it is a $\gamma$-space (Theorem 4 of \cite{GN} and its Corollary), and since the countable union of compact $\gamma$-spaces  is a $\gamma$-space (The union of two compact $\gamma$-spaces is a compact Rothberger space and thus a compact $\gamma$ space. Now apply Jordan's theorem \cite{Jordan} Corollary 14), these topological groups are in fact $\sigma$-compact (thus $\sigma$-scattered) $\gamma$-groups. 
$\Box$

Thus for any cardinal number $\kappa$, the topological group $\reals^{\kappa}$ contains $\sigma$-compact $\gamma$ subgroups of cardinality $\kappa$. It follows for example that the elements with finite support of any power of the integers is a $\sigma$-compact Rothberger bounded topological group. 

\begin{corollary}\label{rothbdedtorothbomega} For each infinite cardinal number $\kappa$ there is a ${\sf T}_0$ topological group $(G,*)$ of cardinality $\kappa$ which is a $\sigma$-compact Rothberger space in all finite powers.
\end{corollary}

\end{document}